\documentclass[11pt]{article}
\usepackage{amsfonts,latexsym,amsmath,amscd,geometry}
\usepackage{amssymb}
\usepackage{xcolor}
\makeatletter \@addtoreset{equation}{section} \makeatother

\newcommand \nc{\newcommand}
\newtheorem{theorem}{Theorem}[section]
\newtheorem{lemma}[theorem]{Lemma}
\newtheorem{proposition}[theorem]{Proposition}

\newtheorem{definition}[theorem]{Definition}

\newtheorem{remark}[theorem]{Remark}

\nc{\ba}{\begin{array}}\nc{\ea}{\end{array}}
\nc{\be}{\begin{eqnarray}}\nc{\ee}{\end{eqnarray}}
\nc{\beq}{\begin{equation}}\nc{\eeq}{\end{equation}}
\nc{\bex}{\begin{eqnarray*}}\nc{\eex}{\end{eqnarray*}}
\nc{\btm}{\begin{theorem}} \nc{\etm}{\end{theorem}}
\nc{\blm}{\begin{lemma}} \nc{\elm}{\end{lemma}}
\nc{\R}{\mathbb{R}}  \nc{\ld}{\lambda}
\nc{\va}{\varphi}
\nc{\ve}{\varepsilon}

\def\P{\mathcal{P}}

\def\pf{\noindent{\bf Proof.\quad}}

\newcommand \qed {\hfill $\Box$}

\def\pft{\noindent{\bf Proof of Theorem 1.1.\quad}}

\begin{document}

\title{Regularity of Structurally Stable Cusp Singularities for Two Families of Quasilinear Wave-type Equations}
\author{Samuel J. Armstrong \footnote{Department of Mathematics, University of Kansas, Lawrence, KS 66045, U.S.A. Email: armstrong@ku.edu}
\quad
Geng Chen
\footnote{Department of Mathematics, University of Kansas, Lawrence, KS 66045, U.S.A. Email: gengchen@ku.edu}
\quad
Tao Huang
\footnote{Department of Mathematics, Wayne State University, Detroit, MI, 48202, U.S.A. Email: taohuang@wayne.edu}
\quad
Yannan Shen
\footnote{Department of Mathematics, University of Kansas, Lawrence, KS 66045, U.S.A. Email: yshen@ku.edu}
}
\maketitle
\begin{abstract}
In this paper, we study two families of quasilinear equations: Hunter-Saxton type and Camassa-Hom type equations, with a paramerter $\lambda\in(0,1)$ whose solutions form cusp singularities. When $\lambda=1$, the first system becomes the scalar conservation law. The main result of this paper is to give regularity of two types of structurally stable singularities: Type I on the singular curve, Type II at the point where cusp singularity forms, for some $\lambda\in(0,1)$. When $\lambda\rightarrow 1$, our result indicts the $C^{1/3}$ regularity at the point where singularity forms, which agrees with the regularity of the generic pre-shock solution. 

 \bigbreak
\noindent


\end{abstract}

\section{Introduction}

Let's start from a general class of nonlinear wave equations parameterized by $\lambda$, which have the form
\beq\label{spe}
	u_{tx}+f'(u)\, u_{xx}+\lambda\, f''(u)\, u_x^2=g(u, u_x)\,.
\eeq
Here $u=u(x,t)$ is a scalar function defined for $(x,t)\in\mathbb{R}\times \mathbb{R}^+$.
More intuitively, equation \eqref{spe} can be formally written as
\begin{equation*}\label{spe2}
 u_{tx}+(f'(u))^{1-\lambda}\, \bigl[(f'(u))^{\lambda}\,u_{x}\bigr]_x=g(u, u_x)\,.
\end{equation*}
Equation \eqref{spe}  includes several important and interesting models when $\lambda$ takes different values.

\begin{itemize}
\item
When $\lambda=1, \ g(u)=0$, \eqref{spe} can be formally written as $u_t +(f(u))_x = 0$, the scalar hyperbolic conservation law. The solution in general forms discontinuities (shock waves). The global wellposedness of small BV solution has been systematically studied, see \cite{Dafermos}. See other models, such as the short pulse equation \cite{SW} with a source term $g(u)$.


\item
When $\lambda=\frac{1}{2}$,
$f'(u)=u,\ g(u)=0$: $u_{tx}+u\, u_{xx}+ \frac{1}{2}u_x^2=0\,$ is the Hunter-Saxton equation, a simplified model from nematic liquid crystals  \cite{BCon,BHR,HS,HZ95a,HZ95b}, where global $H^1$ Solutions were found \cite{BHR,BZZ,CCST,HS,HZ95a,HZ95b}, after the formation of cusp singularity.
\item

When $\lambda=\frac{1}{2a}$, $f'(u)=u^a$ with any positive integer $a$ and $g(u, u_x)=u^{a+1}-P-Q_x$:
\beq\label{E-2}
\begin{array}{lll}u_{t x}+u^a u_{xx}+\frac{1}{2} u^{a-1} u_x^2- u^{a+1}+P+Q_x=0.
 \end{array}
 \eeq
where
$P:=p*[\frac{2a-1}{2}u^{a-1} u_x^2+ u^{a+1}], \ \ Q:=\frac{a -1}{2}\,p*[ u^{a-2}u_x^3],$ \ \ and  \(p = \frac{1}{2}e^{-|x|}\). 
This equation gives Camassa-Holm equation when $a=1$ ($\lambda=\frac{1}{2}$) with $H^1$ solution \cite{BCZ2, BC2}, and Novikov equation when $a=2$ ($\lambda=\frac{1}{4}$) with $W^{1,4}$ solution \cite{CCCS,CCL}.
\end{itemize}



One common feature of these quasilinear equations is the finite time gradient blowup of solutions even with smooth initial data.
 Therefore, it is natural to investigate the global wellposedness and regularity of the general energy conservative solution for equation \eqref{spe} with varying parameter $\lambda\in(0,1)$. 

When $g=0$, the existence of global conservative H\"older continuous weak solutions to \eqref{spe} with exponent $1-\lambda$, when $\lambda\in(0,\frac{1}{3})$ or $\lambda= \frac{1}{2}$, has been shown in \cite{CS}.  Later in \cite{CSZ}, the existence of global conservative H\"older continuous solutions to \eqref{E-2}, i.e. \eqref{spe} with non-local source $g$, has been constructed with exponent $1-\lambda=1-\frac{1}{2a}$, for the integer $a\geq 1$, i.e. when $\lambda=\frac{1}{2N}$ with any natural number $N$. Here the regularity of solution is especially consistent with earlier results on
Camassa-Holm and Novikov equations.



The global existence of solution for \eqref{spe} with $\lambda\in(\frac{1}{3}, \frac{1}{2})\cup (\frac{1}{2},1)$ can be treated using other method \cite{CS,CCSnew}, but unfortunately the regularity of general conservative solution is still open.

In this paper, we consider the generic regularity for  \eqref{spe} with $g=0$ and \eqref{E-2} with varying parameter $\lambda$. More precisely, we prove that the equations \eqref{spe} and \eqref{E-2}  enjoy much better  structurally stable regularities. More precisely, there are two cases of generic singularities: (Type I) on the singularity curve, and (Type II) at the point where an isolated singularity forms, 

The generic singularity is structurally stable, i.e. it keeps its structure under small initial perturbation. The existence of generic singularity for 
Hunter-Saxton, Camassa-Holm and Novikov equation and other equations were proved by \cite{CCCS,CCD, CCST, LZ} using the framework given first in \cite{BressanChen} for the varitional wave equation, by the application of the Thom\rq{}s transversality theorem and the transformation to a semilinear equation. The same method also applies to
\eqref{spe} or \eqref{E-2} with $\lambda\in(0, \frac{1}{3})$ or $\lambda=\frac{1}{2}$. 

There are currently two frameworks which can be used to study the regularity of generic singularities. The first one was established in \cite{BressanHuangYu} for the general solution of variational wave equation, studying the semilinear equation as in \cite{BressanChen}. In this paper, the regularity of three types of singularities, correpsonding to the (Type I) the singular curve, (Type II) the point where cusp singularity forms, (Type III) the intersection point of two cusp singularities from different familes, are provided. Here  the variational wave equation modelling the nematic liquid crystal is another typical model including finite time cusp singularity. See other related works for the global well-posedness of the conservative solution for the variational wave equation \cite{BC2015,BCZ,BZ,HR}. 

The second framework to study the regularity of generic singularity was first established by \cite{BDSV,BSV1,BSV2,BSV3}, when they study the shock formation of compressible Euler equations.  Except showing the $C^{1/3}$ regularity of generic solution right before the shock (the pre-shock solution), this method also provides detail structure of the generic pre-solution solution. This method can apply to solutions in multiple space dimension. Also see many recent advances after this framework.  Espeically, in \cite{KKY}, this framework has been applied to the Camassa holm equation to study the regularity of generic solution right at the point where cusp singularity forms, which is the Type II singularity in \cite{BressanChen} and our current paper. 


In this paper, we use the first framework to study the regularity of generic solutions for \eqref{spe} with $g=0$ and \eqref{E-2} near the singularities. The advantage of this method comparing to the second framework is that it works for the general generic Cauchy initial data, and also provides detailed structure of solutions near the singular curve.

In this paper, we always assume that \(f(u)\) is a \(C^4\) function with \(f'(0) \geq 0\), and that \(f''(u), f'''(u)\), and \(f^{(4)}(u)\) are Lipschitz continuous.


\paragraph{Existence of generic solution}
Using the framework of \cite{BressanChen}, it is not hard to prove the following global existence of generic solutions for \eqref{spe} with $g=0$  and \eqref{E-2} (when \(g= u^{\frac{1}{2\lambda} + 1} - P - Q_x\)). 

For simplicity, in this paper, for equation \eqref{spe} we always assume that 	
\begin{itemize}
				\item[(1)] \(g(u,u_x) = 0\) and \(\lambda = \frac{1}{n} \) where \(n \geq 2\) is an integer, or  
				\item[(2)]  \(g(u,u_x) = u^{\frac{1}{2\lambda} + 1} - P - Q_x\) and \(\lambda = \frac{1}{2a}\) where \(a \geq 1\) is an integer.
		\end{itemize}
		This includes all known related physical models.
\smallskip\\
{\bf Case 1:} When $g=0$, 
we consider the initial boundary value problems for (\ref{spe}) on the region
$(x,t)\in \mathbb{R}^+\times\mathbb{R}^+$ with initial and boundary conditions
		\beq\label{ID1} 
			u(0,t)=0,\qquad u(x,0)=:{u}_{0}(x)\in W_{loc}^{1,\frac{1}{\lambda}}( {\mathbb R}^+)\,,
		\eeq
		and  a compatibility condition
		\beq\label{CC}
			u_0(0)=0 \quad \hbox{and}\quad u'_0(0)=0\,.
		\eeq
Here $W_{loc}^{1,\frac{1}{\lambda}}( {\mathbb R}^+)$ is the Sobolev space with standard notation. 

\begin{theorem}\label{gr_zero} 
	Suppose \(\lambda = \frac{1}{n} \) where \(n \geq 2\) is an integer. For any \(T > 0\) fixed, there exists an open dense set of initial data 
		\[ 
				\mathcal{D} \subset (C^3(\mathbb{R}^+)  \cap  W^{1, \frac{1}{\lambda}}(\mathbb{R}^+) ),
		\] 
		such that for \(u_0 \in \mathcal{D}\), the solution \(u(x,t)\) for \eqref{spe}\eqref{ID1}\eqref{CC} with $g=0$, is \(C^2\) outside finitely many characteristic curves within \(\mathbb{R}^+ \times [0,T]\). 
\end{theorem}
Here for solution $u(x,t)$, we mean the weak solution defined below:
	\begin{definition}{\bf\text (Weak solution)}\label{def}
The function $u(x,t)$, defined for all $(x,t)\in\mathbb{R}^+\times\mathbb{R}^+$,
is a weak solution for \eqref{spe}, if initial and boundary conditions \eqref{ID1} and \eqref{CC} 
are satisfied pointwisely and
	\begin{itemize}
	\item[i.] The equation (\ref{spe}) is satisfied in the weak sense		
	\beq\label{weak}
			\int_{0}^\infty\int_{0}^\infty \Big\{-u_{x}\, 
			(\phi_t+f'(u)\,\phi_x)+(\lambda-1)\, f''(u)u_x^2\,\phi\Big\}\, dx\, dt=0,
		\eeq
		for any test function $\phi\in C_c^1(\mathbb{R}^+\times\mathbb{R}^+)$.
	\item[ii.] For any fixed $t>0$, the function $u(\cdot,t)$ is in $W_{loc}^{1,\frac{1}{\lambda}}( {\mathbb R}^+)$, hence is locally H\"older continuous with exponent $1-\lambda$ by the Sobolev embedding Theorem. 
	\end{itemize}			
\end{definition}
\ 
\bigskip
\\
{\bf  Case 2:} For \eqref{E-2}, i.e. when
\(g= u^{\frac{1}{2\lambda} + 1} - P - Q_x\), we consider the Cauchy problem with the initial data
\beq\label{Ei}u(0,x)=u_0(x)\in  H^1{(\mathbb{R})}\cap W^{1,\frac{1}{\lambda}}{(\mathbb{R})}.\eeq

Similar as in \cite{CCCS} using the frame work of \cite{BressanChen}, we can show the following theorem.
\begin{theorem}\label{gr_nonzero}
 For any \(T > 0\) fixed, there exists an open dense set of initial data 
		\[ 
				\mathcal{D} \subset (C^3(\mathbb{R})  \cap H^1(\mathbb{R}) \cap  W^{1, \frac{1}{\lambda}}(\mathbb{R}) ),
		\] 
		such that for \(u_0 \in \mathcal{D}\), the solution  \(u(x,t)\)  for  \eqref{spe} with \(g= u^{\frac{1}{2\lambda} + 1} - P - Q_x\) is \(C^2\) outside finitely many characteristic curves within \(\mathbb{R} \times [0,T]\). 
\end{theorem}
Here the solution  \(u(x,t)\)   means the weak solution as following
\begin{definition}\label{def1}
Let  $a\in \mathbb{Z}^+$, $\lambda=\frac{1}{2a}$.
The weak solution $u=u(x,t)$ of  the Cauchy problem (\ref{E-2}) (\ref{Ei})  satisfies
\begin{itemize}
\item[(i)] For any fixed $t\geq 0$, $u(\cdot,t)\in H^1{(\mathbb{R})}\cap W^{1,\frac{1}{\lambda}}{(\mathbb{R})}$.
The map $t\mapsto u(\cdot,t)$ is Lipschitz continuous under $L^\frac{1}{\lambda}(\mathbb{R})$ distance.

\item[(ii)] Solution $u=u(x,t)$ satisfies initial condition  (\ref{Ei})  in $L^\frac{1}{\lambda}(\mathbb{R}),$ and
		\begin{equation}\label{nv_weak}
				\begin{split}
						\int_0^\infty \int_\mathbb{R} 
	\left\{-u_x\,  \phi_t -u_x u^{a} \phi_x  +[\frac{1-2a}{2}u^{a-1}u_x^2-u^{a+1} + P+Q_x]  \phi \right\}\, dx\, dt \\  +\int_{\mathbb{R}}(u_{0})_{x}\phi(0,x)\,dx=0 
\end{split}
\end{equation}
for every test function $\phi\in C_c^1(\mathbb{R}\times \mathbb{R}^+)$.

\end{itemize}
\end{definition}

The generic solutions in Theorem \ref{gr_zero} and \ref{gr_nonzero} are all conservative solutions.
The global existence of general energy conservative weak solution can be found in \cite{CS} and \cite{CSZ} respectively. The uniqueness of solution can be proved using the same method of \cite{CCST,BCZ2,CCCS}. We omit it here.

Now we state the main theorem. 

\begin{theorem}\label{main_thm}
		Consider equation \eqref{spe} with generic initial data \(u_0 \in \mathcal{D}\)  given in Theorem \ref{gr_zero} or \ref{gr_nonzero}, and assume additionally that \(u_0 \in C^\infty\).  Let  \(u = u(x,t)\) to denote the solution to \eqref{spe} and let \(\P = (x_0,t_0)\) be a singular point for some $t_0>0$.  
		Then 
		\begin{enumerate}
			\item[(i)] if \(\P\) is a Type I singular point, then there exists a constant \(C_1 = C_1(\lambda) \neq 0\) such that 
					\begin{equation}\label{result_I}
						|u(x,t) - u(x_0,t_0)| =   C_1  |x - x_0 - f'(u(x_0,t_0))(t-t_0)|^{\frac{1}{1+\lambda}} + \mathcal{O}(|x-x_0| + |t-t_0|).
				\end{equation}
			\item[(ii)] If \(\P\) is a Type II singular point, then there exists a constant \(C_2 = C_2(\lambda) \neq 0\) such that 
					\begin{equation} \label{result_II}
							u(x,t) = u(x_0,t_0) + C_2 (x - x_0 - f'(u(x_0,t_0))(t-t_0))^{\frac{2-\lambda}{2+\lambda}} + \mathcal{O}((|x-x_0| + |t-t_0|)^{\frac{2}{2+\lambda}}). 
						\end{equation} 
		\end{enumerate}
\end{theorem}

\begin{remark} We give several remarks.
\begin{itemize}
\item[1.] The Type I (singular curve) and Type II (point where singularity forms) singularities will be precisely defined in Section \ref{genericsingdef}.

\item[2.]  One main observation is that the structurally stable generic regularity is varying with the parameter $\lambda$. It is reasonable to conject that similar conclusions may be valid for all $\lambda\in(0,\frac12].$

\item[3.] When $\lambda\rightarrow 1$, the Type II generic regularity, corresponding to the point when cusp singularity forms, the Hölder exponent approaches $\frac{1}{3}$, which agrees with the result for the pre-shock of Euler equations using the frameworks in \cite{BDSV,BSV1,BSV2,BSV3}. 
When $\lambda=\frac{1}{2}$, the regularity for Type II regularity agrees with the result for \cite{KKY} for the Camassa-Holm equation with speical initial data. 

\end{itemize}

\end{remark}

The remainder of the paper can be divided into three main sections. In Section 2, we review the main equations  and ideas from \cite{CS} while establishing the existence of solutions of \eqref{spe} for generic smooth initial data.
In Section 3, we prove Theorem \ref{main_thm} when \(g(u,u_x) = 0\) by introducing a new change of variables and computing Taylor expansions of the generic regular solutions about the singular points. Finally, in Section 4, we revisit Theorem \ref{main_thm} for the case where \(g(u,u_x) \neq 0\) and look at the special cases for the Camassa-Holm and Novikov equations.

%

\section{General type Hunter-Saxton type equation with $g(u,u_x)=0$}

In this section, we first investigate the Cauchy problem for the general Hunter-Saxton type equation \eqref{spe} with $g(u,u_x)=0$
\beq\label{LGHS}
u_{tx}+f'(u)u_{xx}+\lambda f^{''}(u)u_x^2=0,
\eeq
equipped with initial data \(u_0(x) \in W^{1,\frac{1}{\lambda}}(\mathbb{R}^+)\). 

Multiplying the equation \eqref{LGHS} by $u_x^{\frac{1}{\lambda}-1}$ and integrating over $\R$, we have
\beq\notag
\lambda\frac{d}{dt}\int (u_x)^{\frac{1}{\lambda}}\,dx
+\int f'(u)u_{xx}(u_x)^{\frac{1}{\lambda}-1}\,dx
+\int \lambda f''(u)(u_x)^{\frac{1}{\lambda}+1}=0.
\eeq
Integrating by parts implies
\beq\notag
\int f'(u)u_{xx}(u_x)^{\frac{1}{\lambda}-1}\,dx
=\lambda\int f'(u)\frac{\partial}{\partial x}\left(u_x^{\frac{1}{\lambda}}\right)
=-\lambda\int f''(u)(u_x)^{\frac{1}{\lambda}+1}.
\eeq
Combining these two estimates, it holds
\beq\notag
\lambda\frac{d}{dt}\int (u_x)^{\frac{1}{\lambda}}\,dx=0
\eeq
which provides the energy conservation for any regular solution of \eqref{LGHS}.

Note, $ (u_x)^{\frac{1}{\lambda}}$ and  $(u_x)^{\frac{1}{\lambda}-1}$ make sense as  $ (u_x^2)^{\frac{1}{2\lambda}}$ and  $(u_x^2)^{\frac{1}{2\lambda}-\frac{1}{2}}$.

\subsection{Semi-linear system}
The equation of characteristic of \eqref{LGHS} is 
\beq\notag
\frac{d x^c(s)}{ds}=f'(u(x^c(s),s)).
\eeq
We denote the characteristic through the point $(x,t)$ as $x^c(s,x,t)$, or equivalently \(t^c(r,x,t)\),  and introduce a change of coordinates $(x,t)\rightarrow(Y,\tau)$ as 
\[
Y=\begin{cases}
		\int_0^{x^c(0,x,t)}\left(1+u_x^2(y,0)\right)^{\frac{1}{2\lambda}}\,dy, & \text{when the characteristic passing \((x,t)\) intersects \(t = 0\)},  \\
		-t^c(0,x,t)f'(0), & \text{when the characteristic passing \((x,t) \) intersects \(x = 0\). }
	\end{cases}
\]
and 
\[ 
	\tau = t. 
\] 
Direct computation implies 
\beq\notag
Y_t+f'(u)Y_x=0,\quad \tau_t=1,\quad \tau_x=0,
\eeq
from which,  it holds for any smooth function $m$
\beq\label{devcomp}
m_t+f'(u)m_x=m_{\tau},\quad m_x=m_YY_x.
\eeq

To rewrite the nonlinear equation \eqref{LGHS} into a semi-linear system, we introduce the following dependent variables
\beq\label{defVX}
v=2\arctan u_x,\quad \xi=\frac{(1+u_x^2)^{\frac{1}{2\lambda}}}{Y_x}.
\eeq
By direct computations, we obtain the following semi-linear system, which is equivalent to \eqref{LGHS}
\beq\label{semieq}
\begin{cases}
		u_Y=\frac{1}{2}\xi\sin v\left(\cos^2\frac{v}{2}\right)^{\frac{1}{2\lambda}-1},\\
v_{\tau}=-2\lambda f''(u)\sin^2\frac{v}{2},\\
\xi_{\tau}=\frac12\xi f''(u)\sin v.
\end{cases}
\eeq
We impose the initial data 
\begin{equation}\label{id_zero}
	\begin{cases}
		u(Y,0) := u_0(x(Y,0)), \\
		v(Y,0) := 2\arctan (u'_0(x(Y,0))),\\
		\xi(Y,0) := 1,
	\end{cases}
\end{equation}
and the boundary conditions 
\begin{equation}\label{bc_zero}
		\begin{cases}
	u(-f'(0)\tau,\tau) := 0, \\
	v(-f'(0)\tau,\tau) := 0, \\
	\xi(-f'(0)\tau,\tau) := 1.
\end{cases}
\end{equation}
By \eqref{devcomp},\eqref{defVX} and \eqref{semieq}, we obtain
\beq\label{xtY}
\begin{cases}
x_Y=\xi\left(\cos^2\frac{v}{2}\right)^{\frac{1}{2\lambda}},\\
x_\tau=f'(u),
\end{cases}
\quad
\begin{cases}
t_Y=0,\\
t_{\tau}=1.
\end{cases}
\eeq

\begin{remark}
(1) The system \eqref{semieq} is invariant under the translation by $2\pi$ in $v$. Hence, we consider $v\in[-\pi,\pi]$ for simplicity. 

(2) When $\frac{1}{\lambda}$ is not a positive integer, the right side of the first equation in \eqref{semieq} is not smooth since $\left(\cos\frac{v}{2}\right)^{\frac{1}{\lambda}-1}\in C^{k_0,\alpha}$, where $k_0=\left[\frac{1}{\lambda}-1\right]$ and $\alpha=\frac{1}{\lambda}-1-\left[\frac{1}{\lambda}-1\right]$.


\end{remark}
By assuming \(v \in [-\pi, \pi]\), it is convenient to write 
\begin{equation}\label{uandx}
		u_Y = \xi \sin \frac{v}{2} \left(\cos \frac{v}{2} \right)^{\frac{1}{\lambda} - 1} \quad \text{ and } \quad x_Y = \xi \left( \cos \frac{v}{2} \right)^{\frac{1}{\lambda}} 
\end{equation}
since \(\cos \frac{v}{2} \geq 0\). For better regularity, we introduce another dependent variable 
\beq\notag
S=u_t+f'(u)u_x=u_\tau.
\eeq
Direct computation implies
\beq\notag
\begin{split}
&S_Y=\frac{u_{\tau x}}{Y_x}\\
=&\frac{1}{Y_x}\left(u_{tx}+f'(u)u_{xx}+f''(u)u_x^2\right)\\
=&\frac{(1-\lambda)}{Y_x}f''(u)u_x^2\\
=&(1-\lambda)\xi f''(u)\sin^2\frac{v}{2}\left(\cos\frac{v}{2}\right)^{\frac{1}{\lambda}-2}
\end{split}
\eeq
Thus the system \eqref{semieq} can also be written as follows
\beq\label{semieqS}
\begin{cases}
u_\tau=S\\
S_Y=(1-\lambda)\xi f''(u)\sin^2\frac{v}{2}\left(\cos\frac{v}{2}\right)^{\frac{1}{\lambda}-2},\\
v_{\tau}=-2\lambda f''(u)\sin^2\frac{v}{2},\\
\xi_{\tau}=\frac12\xi f''(u)\sin v.
\end{cases}
\eeq

\bigskip
When $\lambda \in (0,\frac{1}{3}] \cup \{\frac{1}{2}\}$, the right-hand side of \eqref{semieqS} is Lipschitz continuous. When \(\lambda \in (0,\frac{1}{2}]\), the right-hand side of \eqref{semieq} is Lipschitz continuous. 

\bigskip

\subsection{Generic Regularity}
The introduction of \(v = 2\arctan u_x\) captures the points where \(|u_x| \to \infty\) along the level sets \(\{v = \pi\}\) and \(\{v = -\pi\}\). Moreover, the Jacobian 
\begin{equation}\label{jacobian}
	\begin{vmatrix}
		x_Y & t_Y \\
		x_\tau & t_\tau 
\end{vmatrix} = \xi (\cos \frac{v}{2})^{\frac{1}{\lambda}}
\end{equation}
is nonzero when \(v \neq \pm \pi\). That is, solutions \(u\) will preserve smoothness in neighborhoods of points \((x,t) = (x(Y,\tau), t(Y,\tau))\) such that \(v(Y,\tau) \neq \pm \pi \). Without loss of generality, it suffices to study the level set \( \{v = \pi\}\) to investigate the singularities of solutions. To this end, Theorem \ref{gr_zero} states that piecewise \(C^2\) solutions of \eqref{spe} are dense for smooth enough initial data, and that the level set \(\{v=\pi\}\) for such solutions are made up of finitely many smooth curves.

Now we come to prove Theorem \ref{gr_zero}. We first present two important lemmas related to the peturbations of solutions and the smoothness of solutions on compact domains.  
\begin{lemma}\label{perturb}
		Let \((u,v,\xi)\) be a smooth solution of \eqref{semieq} and \( (Y_0,\tau_0) \in \mathbb{R^+} \times \mathbb{R}^+\).
\begin{enumerate}
\item[(1)] If \((v,v_Y,v_{YY})(Y_0,\tau_0) = (\pi,0,0)\), then there exists a 3-parameter family of smooth solutions \( (u^\theta, v^\theta, \xi^\theta)\) of \eqref{semieq} depending smoothly on \(\theta \in \mathbb{R}^3\) with \(\theta\) in a small enough open ball centered at \((0,0,0)\) such that 
\begin{enumerate}
	\item[(i)] when \(\theta = 0\), one recovers the original solution \((u,v,\xi)\), and
	\item[(ii)] at the point \( (Y_0,\tau_0)\), when \(\theta = 0\) one has \(\textbf{rank }D_\theta (v^\theta, v_Y^\theta, v_{YY}^\theta) = 3\). 
\end{enumerate}
\item[(2)] If \((v,v_Y,v_\tau)(Y_0,\tau_0) = (\pi,0,0) \), then there exists aa 3-parameter family of smooth solutions \((u^\theta, v^\theta, \xi^\theta)\) of \eqref{semieq} depending smoothly on \(\theta \in \mathbb{R}^3\) satisfying (i)-(ii) above, and at \((Y_0,\tau_0)\) when  \(\theta = 0\) one has \( \textbf{rank }D_\theta (v^\theta, v_Y^\theta, v_{\tau}^\theta) = 3 \).
\end{enumerate}

\end{lemma}

\pf
Let \((u,v,\xi)\) be a smooth solution to the semilinear system. We construct families of solutions \((u^\theta, v^\theta, \xi^\theta)\) to the semilinear system with perturbations on the initial data as 
\[ 
		\begin{cases}
				u^\theta(Y,0) & = u(Y,0) + \sum_{i=1}^3 \theta_i U_i(Y,0)\\
				v^\theta(Y,0) & = v(Y,0) + \sum_{i=1}^3 \theta_i V_i(Y,0)\\
				\xi^\theta(Y,0) & = \xi(Y,0) + \sum_{i=1}^3  \theta_i \Xi_i(Y,0)\\
		\end{cases}	
	\] 
	respectively for some suitable functions \(U_i(Y,0) ,V_i(Y,0), \Xi_i(Y,0) \in C_c^\infty (\mathbb{R}^+) \), subject to the boundary conditions 
	\[ 
			\begin{cases}
	u^\theta(-f'(0)\tau,\tau) := 0, \\
	v^\theta(-f'(0)\tau,\tau) := 0, \\
	\xi^\theta(-f'(0)\tau,\tau) := 1.
\end{cases}
	\] The perturbed functions must satisfy 
	\[
		\begin{cases}
				u_{\tau}^\theta =  S & =: f_1^\theta, \\
				v_\tau^\theta  = -2 \lambda f''(u^{\theta}) \sin^2 \frac{v^\theta}{2} & =: f_2^\theta, \\ 
				\xi_\tau^\theta  = \frac{1}{2}f''(u^{\theta}) \xi^\theta \sin v^\theta & =: f_3^\theta.
		\end{cases}
	\] 	
	Taking derivatives of \(v_\tau^\theta\) with respect to \(Y\) yields 
\[ 
			v^{\theta}_{Y \tau }   = -\lambda\left(2 f'''(u^{\theta}) \xi^\theta \sin^3 \frac{v^\theta}{2} (\cos \frac{v^\theta}{2})^{\frac{1}{\lambda} -1}  + f''(u^\theta)  \sin v^\theta  v^\theta_Y\right)  := f_4^\theta  \\
\] 
and
\begin{multline*} 
		v^\theta_{YY \tau}  = - \lambda \Biggl( 2f^{(4)}(u^\theta) (\xi^\theta)^2 \sin^4 \frac{v^\theta}{2} (\cos \frac{v^\theta}{2})^{\frac{2}{\lambda} - 2} \\
		  + 2f'''(u^\theta) \xi_Y^\theta  \sin^3     \frac{v^\theta}{2} (\cos \frac{v^\theta}{2})^{\frac{1}{\lambda} -1} + 3 f'''(u^\theta) \xi^\theta \sin^2 \frac{v^\theta}{2} (\cos \frac{v^\theta}{2})^{\frac{1}{\lambda}}v_{Y}^\theta \\
		  - \frac{1-\lambda}{\lambda} f'''(u^\theta) \xi^\theta \sin^4 \frac{v^{\theta}}{2} (\cos \frac{v^\theta}{2})^{\frac{1}{\lambda} - 2}v_{Y}^\theta + f'''(u^\theta) \xi^\theta \sin \frac{v^\theta }{2} (\cos \frac{v^\theta}{2})^{\frac{1}{\lambda} -1}\sin v^\theta v_Y^\theta  \\
		 + f''(u^\theta) \cos v^\theta (v_Y^\theta)^2 + f''(u^\theta) \sin v^\theta v_{YY}^\theta   \Biggr) =: f_5^\theta . 
\end{multline*} 
Lastly,
\[ 
		\xi^\theta_{Y\tau} = \frac{1}{2} \left( f'''(u^\theta) (\xi^\theta)^2 \sin^3 \frac{v^\theta}{2} (\cos \frac{v^\theta}{2})^{\frac{1}{\lambda} - 1} \sin v^\theta  + f''(u^\theta) \xi_Y^\theta \sin v^\theta +  f''(u^\theta) \xi^\theta \cos v^\theta v^\theta_Y \right) =: f^\theta_6.  
\] 
This gives us the complete ODE system 
\[ 
	\frac{\partial}{\partial \tau} \begin{pmatrix}
		u^\theta \\
		v^\theta \\
		\xi^\theta  \\
		v_Y^\theta \\
		v_{YY}^\theta\\
		\xi_Y^\theta \\
		\end{pmatrix} = \begin{pmatrix}
				f_1^\theta \\
				f_2^\theta \\
				f_3^\theta \\
				f_4^\theta \\
				f_5^\theta \\
				f_6^\theta 
	\end{pmatrix}.
\] 
When \(\lambda \in (0,\frac{1}{3}] \cup \{\frac{1}{2}\} \) the right-hand side is Lipschitz continuous, so we can choose suitable peturbations \(V_i\) for which \(\textbf{rank } D_\theta(v^\theta, v_Y^\theta, v_{YY}^\theta) = 3 \) at \(\theta = 0\) when \((Y,\tau) = (Y_0,\tau_0)\) (see Lemma 4 from \cite{BressanChen} for details).  

On the other hand, we can also compute 
\[ 
	v^\theta_{\tau \tau} = -2\lambda Sf'''(u^\theta) \sin^2 \frac{v^\theta}{2} + 4\lambda^2 (f''(u^\theta))^2 \sin^3 \frac{v^\theta}{2} \cos v^\theta  =: f^\theta_7. 
\] 
This gives the complete ODE system 
\[ 
	\frac{\partial }{\partial \tau} \begin{pmatrix}
			u^\theta \\
		v^\theta \\
		\xi^\theta  \\
		v_Y^\theta \\
		v_{\tau}^\theta\\
		\end{pmatrix} = \begin{pmatrix}
				f_1^\theta \\
				f_2^\theta \\
				f_3^\theta \\
				f_4^\theta \\
				f_7^\theta \\
\end{pmatrix}.
\] 
Similar to above, the right-hand side is Lipschitz continuous when \(\lambda \in (0,\frac{1}{3}] \cup \{\frac{1}{2}\}\), so we can choose suitable peturbations for which  \(\textbf{rank } D_\theta(v^\theta, v_Y^\theta, v_{\tau}^\theta) = 3 \) when \(\theta = 0\) and \((Y,\tau) = (Y_0,\tau_0)\).
\qed

\begin{lemma} \label{density}
	For any \(T,M > 0\), let  
	\[ 
		\Omega = \{  (Y,\tau):  -f'(0)\tau \leq Y \leq  M, \; 0 \leq \tau \leq T\}
	\] 
	and define \(\mathcal{S}\) to be the family of all \(C^2\) solutions \((u,v,\xi)\) to \eqref{semieq}. Define \(\mathcal{S}' \subset \mathcal{S}\) to be the subfamily of solutions \((u,v,\xi)\) such that for any \((Y,\tau) \in \Omega\), the values
	\begin{equation} \label{values}
			(v,v_{Y},v_{YY}) = (\pi,0,0), \quad (v,v_Y,v_{\tau}) = (\pi,0,0) 
\end{equation} 
cannot be attained. Then \(\mathcal{S}'\) is a relatively open and dense subset of \(\mathcal{S}\), in the topology induced by \(C^2(\Omega)\). 
\end{lemma}
We omit the proof of Lemma \ref{density} for brevity, as it follows in the same manner as the proofs presented in \cite{BressanChen,LZ}. With Lemma \ref{perturb} and Lemma \ref{density}, we can now prove Theorem \ref{gr_zero}.

\pft 
Define the space
\[ 
		\mathcal{M} := C^3 (\mathbb{R}^+) \cap W^{1,\frac{1}{\lambda}}(\mathbb{R}^+)
\] 
equipped with the norm 
\[ 
		\|u\|_{\mathcal{M}} := \|u\|_{C^3} + \|u\|_{W^{1,\frac{1}{\lambda}}}.
\] 
Pick any initial data \(\hat{u}_0 \in \mathcal{M}\) and for \(\delta > 0\), denote the open ball
\[ 
		B_\delta := \{ u_0 \in \mathcal{M}: \|u_0 - \hat{u}_0\|_{\mathcal{M}} < \delta\}.
\] 
For \(u_0 \in \mathcal{M}\), we know that 
\[ 
	u_0(x) \to 0, \quad u_{0,x}(x) \to 0
\] 
as \(x \to \infty\). Specifically, there exists some \(\rho > 0\) large enough so that \(u_0\) and \(u_{0,x}\) remain uniformly bounded for all \(x \geq \rho\). By a standard comparison argument, the corresponding solution \(u\) must have \(u_x\) remain uniformly bounded on \(\{ (x,t): t \in [0,T], \; x \geq \rho \}\). In particular, \(u\) remains \(C^2\) outside of the set
\[
		\mathcal{N} := [0, \rho] \times [0,T]. 
\]
Thus, the singularities of \(u\) can only form within \(\mathcal{N}\). Now define the map \(\Lambda\) by
\[
		(Y,\tau) \mapsto \Lambda(Y,\tau) :=  (x(Y,\tau), t(Y,\tau)).
\]
Using the set \(\Omega\) defined in Lemma \ref{density}, choose \(M  > 0\) sufficiently large and   \( \delta > 0\) sufficiently small  so that for all \(u_0 \in B_\delta\), we have \( \mathcal{N} \subset \Lambda(\Omega)\). Finally, we will define \(\mathcal{D}\) as follows: We say \(u_0 \in \mathcal{D}\) if \(u_0 \in B_\delta\) and the corresponding solution \((u,v,\xi)\) to \eqref{semieq}-\eqref{bc_zero} does not attain the values \eqref{values} for any \((Y,\tau)\) such that \(\Lambda(Y,\tau) \in \mathcal{N}\).

We first show \(\mathcal{D}\) is open in \(C^3\). Let \( u_0^\nu \in B_\delta\)  be a sequence of functions with \(\nu \geq 1\) that converges to \(u_0\) with \(u_0^\nu \not \in \mathcal{D}\). Then there exists a seuqence of points \((Y^\nu, \tau^\nu)\) such that the corresponding solution \( (u^\nu, v^\nu, \xi^\nu)\) satisfy 
\[ 
		(v^\nu, v^{\nu}_Y, v^\nu_{YY}) = (\pi,0,0), \quad (x(Y^\nu,\tau^\nu), t(Y^\nu,\tau^\nu)) \in \mathcal{N}. 
\] 
for all \(\nu \geq 1\). Note that \(\mathcal{N}\) is compact, so there exists a subsequence of \( (Y^{\nu_k}, \tau^{\nu_k})  \) of \( (Y^{\nu}, \tau^{\nu}) \) and a point \((\overline{Y}, \overline{\tau} ) \) for which \( (Y^{\nu_k}, \tau^{\nu_k})  \to (\overline{Y}, \overline{\tau} ) \). But by continuity, the corresponding solution \((u,v,\xi)\) from the initial data \(u_0\) must satisfy 
\[ 
		(v,v_Y,v_{YY})(\overline{Y},\overline{\tau}) = (\pi,0,0), \quad (x(\overline{Y},\overline{\tau}), t(\overline{Y},\overline{\tau})) \in \mathcal{N}. 	
\] 
Thus, \(u_0 \not \in \mathcal{D}\). Repeating the same argument when \((v,v_{Y},v_\tau) = (\pi,0,0)\)  shows that \(\mathcal{D}\) is open. 

Next, we show that \(\mathcal{D}\) is dense. By a small perturbation, assume that \(u_0 \in C^\infty\). By Lemma \ref{density}, construct a sequence of solutions \( (u^\nu, v^\nu, \xi^\nu)\) for \(\nu \geq 1\) where 
\begin{enumerate}
		\item[(i)] the values \eqref{values} are never attained for any \((Y,\tau) \in \Omega\).
			\item[(ii)] The \(C^k\) norm satisfies 
			\[ 
					\lim_{\nu \to \infty} \|(u^\nu - u, v^\nu - v, \xi^\nu - \xi) 	\|_{C^k(I)} = 0
			\] 
			for all \(k \geq 1\) and all bounded sets \(I \subset \mathbb{R}^+ \times [0,T]\). 
\end{enumerate}	
In particular, at \(t = 0\), we have 
\[
	\lim_{\nu \to \infty} \|u_0^\nu - u_0\|_{C^k([a,b])} = 0
\] 
for all bounded intervals \([a,b] \subset \mathbb{R}^+\). To get convergence in the far field, introduce the cutoff function \(\eta \in C^\infty_c \) where 
\[ 
	\eta(x) = \begin{cases}
			1, & \text{ if } 0 \leq x \leq \ell, \\
			0, & \text{ if } x \geq \ell + 1,
	\end{cases}
\] 
for \(\ell > 0\). For each \(\nu \geq 1\), consider the new initial data 
\[ 
	\tilde{u}_0^\nu := \eta u_0^\nu + (1-\eta)u_0.
\] 
Then 
\[ 
	\lim_{\nu \to \infty} \|\tilde{u}_0^\nu - u_0\|_{\mathcal{M}} = 0. 
\] 
Now choose \(\ell > \rho \) large enough so that for any \((x,t) \in \mathcal{N}\), 
\[ 
	\tilde{u}^\nu(x,t) = u^\nu(x,t). 
\] 
Then \(\tilde{u}^\nu\) is \(C^2\) outside of \(\mathcal{N}\) and thus \(\tilde{u}_0^\nu \in \mathcal{D}\) for large enough \(\nu\). We conclude \(\mathcal{D}\) is dense in \(B_\delta\). 

Finally, let \(u_0 \in \mathcal{D}\). We wish to show the corresponding solution \(u\) is piecewise \(C^2\) within \(\mathcal{N}\). Let \((Y_0, \tau_0) \in \Omega \). If \(v(Y_0,\tau_0) \neq \pi\), then the coordinate change \( (x,t) \mapsto (Y,\tau)\) is locally invertible by \eqref{jacobian}. Thus, suppose that \(v(Y_0,\tau_0) = \pi\). By definition of \(u_0 \in \mathcal{D}\), we must have \(v_\tau(Y_0,\tau_0) \neq 0\) or \(v_Y(Y_0,\tau_0) \neq 0\). By continuity, there exists some \(\epsilon  > 0\) such that the values \eqref{values} are not attained in the open neighborhood 
\[ 
		\Omega' := \{(Y,\tau): -f'(0) \tau - \epsilon < Y < M +\epsilon, -\epsilon < t < T+\epsilon \}.
\] 
By the implicit function theorem, the set 
\[
		S^v : = \{(Y,\tau) \in \Omega': v(Y,\tau) = \pi\} 
\] 
is a one-dimensional embedded mainfold of class \(C^2\). We claim that the number of connected components of \(S^v\) that intersect \(\Omega\) is finite. Suppose, to the contrary, that \(P_n \in S^v \cap \Omega\) is a sequence of points for \(n \geq 1\) belonging to distinct components. Choose a subsequence \(P_{n_k}\) of \(P_n\) such that \(P_{n_k} \to \overline{P}\) for some \(\overline{P} \in S^v \cap \Omega\). By assumption, \((v_Y,v_\tau)(\overline{P}) \neq (0,0)\). Thus, by the implicit function theorem, there exists a neighborhood \(U\) of \(\overline{P}\) such that \(\gamma := S^v \cap U\)  is a connected \(C^2\) curve. That is, \(P_{n_k} \in \gamma\) for large enough \(k\), contradicting \(P_{n_k}\) belonging to distinct components.   

From our previous argument, there are only finitely many points in \(\Omega\)  where \(v = \pi\), \(v_Y \neq 0\), and \(v_\tau = 0\). Call these points \(P_i = (Y_i,\tau_i)\) for \(i = 1,\hdots, m\). Similarly, there are finitely many points in \(\Omega\) such that \(v = \pi\), \(v_{Y} = 0\), and \(v_{\tau} \neq 0\). Call these points \(Q_j = (Y_j',\tau_j')\) where \(j = 1,\hdots, n\). Then \(S^v \backslash \{P_1,\hdots, P_m,Q_1,\hdots, Q_n\}\) has finitely many connected components that intersect \(\Omega\). Any one of these components \(\gamma_k\) is a connected curve such that \(v = \pi, v_Y \neq 0\) for any \((Y,\tau) \in \gamma_k \). Thus, we can write 
\[ 
	\gamma_k = \{(Y,\tau): Y = \phi_k(\tau), a_k < \tau < b_k\} 
\] 
for a suitable \(C^1\) function \(\phi_k\). As a result, the image \(\Lambda(\gamma_k)\) is a \(C^2\) curve in the \(x\)-\(t\) plane. Moreover, the singular set \(\Omega(S^v)\) is the union of finitely many points \(p_i = \Lambda(P_i)\) and \(q_j = \Lambda(Q_j)\) together with finitely many \(C^2\) curves \(\Lambda(\gamma_k)\). This proves Theorem \ref{gr_zero}.
\qed

\subsection{Generic singularities}
\label{genericsingdef}
From Theorem \ref{gr_zero}, generic regular solutions can have two types of behavior near singular points $(Y_0,\tau_0)$ where \(v(Y_0,\tau_0) = \pi\): 
\beq\label{tp1}
\mbox{Type I: }\quad  v=\pi,\quad  v_Y\neq 0,
\eeq
\beq\label{tp2}
\mbox{Type II: }\quad v=\pi,\quad  v_Y= 0, \quad v_{YY}\neq 0.
\eeq

If we assume $f''(u)\neq 0$ for any $u\in\R$, then it holds $v_{\tau}\neq 0$ by the second equation in system \eqref{semieq}.

Points of Type 1 form a locally finite family of $C^2$ curves in the $Y$-$\tau$ plane. Their images form a family of characteristic curves in $x$-$t$ plane where the solution $u = u(x, t)$ is singular i.e., $u$ is not diﬀerentiable.
Points of Type 2 are isolated. Their images in $x$-$t$ plane are points where two singular curves initiate or terminate.


\section{Regularity Near Singularities when \(g(u,u_x) = 0\)} 
We now prove Theorem \ref{main_thm} for the case where \(g(u,u_x) = 0\). When \(\lambda = \frac{1}{n}\) where \(n \geq 2\) is an integer, solutions \( (u,v,\xi)\) to \eqref{semieq} with smooth initial data will remain smooth in the \(Y\)-\(\tau\) plane. As introduced in \cite{BressanHuangYu}, the local behavior of the solution \(u\)  near the singular points can be computed by the Taylor expanding \(u\) and \(x\) in the \(Y\)-\(\tau\) variables. We will demonstrate a different, but similar,  approach for \(\lambda = \frac{1}{n}\) where \(n \geq 2\) is an integer. Rather than computing all derivatives explicitly, we introduce a change of variables that allows us to compute just one spatial derivative for both \(u\) and \(x\). Going forward, we will denote \(\P = (Y_0,\tau_0)\) as a singular point where \(v(Y_0,\tau_0) = \pi\) and let \(x_0 = x(Y_0,\tau_0)\), \(t_0 = t(Y_0,\tau_0) \), \(\xi_0 = \xi(Y_0,\tau_0)\). Moreover, we will assume \(u_0 \in \mathcal{D} \cap C^\infty \), so we may take as many derivatives as we need on \(u, v,\) and \(\xi\). 

\subsection{Type I Singularity}
Suppose \(\P\) is a Type I singular point. That is, \(v_Y(Y_0,\tau_0) \neq 0\). Introduce the change of variables 
\begin{equation} 
		W = (Y-Y_0)^{\frac{1}{\lambda}}, \quad Z = (Y-Y_0)^{\frac{1 + \lambda }{\lambda  }}.
\end{equation} 
In particular, \((W,\tau) = (0,\tau_0)\) when \((Y,\tau) = (Y_0, \tau_0)\). It is important to note that when \(\lambda = \frac{1}{n}\), we have \(\frac{1}{\lambda} = n\) and \(\frac{1 + \lambda}{\lambda} = n+1\). Thus, the powers are of opposite parity. First, suppose \(n\) is even.  When \(Y-Y_0 > 0\), we have 
\[
		Y_W = \frac{\lambda}{ W^{1-\lambda}}, \quad Y_Z = \frac{\lambda}{(1+\lambda) Z^{1 - \frac{\lambda}{1 + \lambda} }}
\] 
away from \(W = 0\) and \(Z = 0\). 
We can use the chain rule to compute 
\begin{align*}
		u_W(W,\tau) & = u_Y(W, \tau) Y_W \\
					& = \lambda \xi(W,\tau) \sin \frac{v(W, \tau)}{2} \frac{(\cos \frac{v(W, \tau)}{2} )^{\frac{1}{\lambda} - 1}}{W^{1-\lambda}}\\
					& = \lambda \xi(W,\tau) \sin \frac{v(W, \tau)}{2}  \left( \frac{\cos \frac{v(    W, \tau)}{2} }{W^{\lambda} }\right)^{\frac{1}{\lambda} - 1}
\end{align*}
and 
\begin{align*}
		x_Z(Z,\tau) & = x_Y(Z,\tau) Y_Z \\
				 & = \frac{\lambda}{1 + \lambda} \xi(Z,\tau) \frac{(\cos \frac{v(Z,  \tau)}{2} )^{\frac{1}{\lambda}} }{Z^{1 - \frac{\lambda}{1 + \lambda}}}\\
				 & =  \frac{\lambda}{1 + \lambda} \xi(Z,\tau)  \left(\frac{\cos \frac{v(Z,\tau) }{2}}{Z^{\frac{\lambda}{1+\lambda}}} \right)^{\frac{1}{\lambda}}.
\end{align*}
Here we can see \(u(W,\tau)\) and \(x(Z,\tau)\) are continuously differentiable up to \(W=0\) and \( Z = 0 \) respectively.  At \(\tau = \tau_0\), we can use L'Hopital's rule to compute 
\begin{align*}
		u_W(0,\tau_0) & := \lim_{W \to 0^+} u_W(W,\tau_0) \\
					  & = \lambda \xi_0 \left( \lim_{W \to 0^+} \frac{\cos \frac{v(W, \tau_0)}{2} }{W^{\lambda} }\right)^{\frac{1}{\lambda} - 1}\\
					  & = \lambda \xi_0 \left( \lim_{W \to 0^+}\frac{- \frac{\lambda}{2}v_Y(    W, \tau_0) (\sin \frac{v(W,\tau_0)}{2}) W^{\lambda - 1} }{\lambda W^{\lambda - 1}} \right)^{\frac{1}{\lambda} - 1}\\
					  & = (-1)^{\frac{1}{\lambda} - 1} \frac{\lambda \xi_0 }{2^{\frac{1}{\lambda} - 1} }v^{\frac{1}{\lambda} -1 }_{Y}(Y_0,\tau_0) \neq 0
\end{align*}
and 
\begin{align*}
	x_Z(0,\tau_0) & := \lim_{Z \to 0} x_Z(Z,\tau_0) \\
				& = \frac{\lambda \xi_0}{1+\lambda}\left(\lim_{Z \to 0} \frac{\cos \frac{v(Z,\tau_0) }{2}}{Z^{\frac{\lambda}{1+\lambda}}}  \right)^{\frac{1}{\lambda}}\\
				& = (-1)^{\frac{1}{\lambda}} \frac{\lambda \xi_0}{2^{\frac{1}{\lambda}}(1 + \lambda)}v^{\frac{1}{\lambda}}_Y(Y_0,\tau_0)) \neq 0
\end{align*}
We claim the following: 
\begin{proposition} \label{holder_I}
		Suppose \(\lambda = \frac{1}{n}\) for some even integer \(n \geq 0\). Then there is some \(\delta > 0\) such that \(u_W(\cdot, \tau_0) \in C_W^{0,\lambda}([0,\delta]) \) and \(x_Z(\cdot, \tau_0) \in C_Z^{0, \frac{\lambda}{1+\lambda}}([-\delta,\delta])\). 
\end{proposition}
\pf
		Consider the function 
		\[ 
			h(Y) : = \frac{\cos \frac{v(Y,\tau_0)}{2} }{Y-Y_0}.
		\] 
		It is clear that \(h\) is differentiable up to \(Y = Y_0\) with 
		\[ 
			h'(Y) = \frac{-\frac{1}{2}(Y-Y_0) \sin (\frac{v(Y,\tau_0)}{2}) v_Y(Y,\tau_0) - \cos \frac{v(Y,\tau_0)}{2} }{(Y-Y_0)^2 }.
		\] 
		But we can compute 
		\begin{align*}
				& \lim_{Y \to Y_0} h'(Y) = \\
				& \lim_{Y \to Y_0} \frac{ - \frac{1}{4}(Y-Y_0) \cos(\frac{v(Y,\tau_0)}{2}) v_{Y}^2(Y,\tau_0)   - \frac{1}{2}(Y-Y_0) \sin(\frac{v(Y,\tau_0)}{2})    v_{YY}(Y_0,\tau_0) }{2(Y-Y_0)}\\
				& = -\frac{1}{4}v_{YY}(Y_0,\tau_0).
		\end{align*}
		In particular, \(h\) must be Lipschitz near \(Y_0\). This is noteworthy since 
		\[ 
				u_W(W,\tau_0) = \lambda \xi(W,\tau_0) \sin \frac{v(W,\tau_0)}{2} (h(W^\lambda + Y_0))^{\frac{1}{\lambda}- 1}
		\] 
		and 
		\[ 
				x_Z(Z,\tau_0) = \frac{\lambda}{1+\lambda} \xi(Z,\tau_0) (h(Z^{\frac{\lambda}{1+\lambda}}+ Y_0))^{\frac{1}{\lambda}}. 
		\] 
		From here, we can see  \(h^{\frac{1}{\lambda}} \in C_Z^{0,\frac{\lambda}{1+\lambda}}([-\delta,\delta]) \) for some \(\delta > 0\). If \(\lambda \in (0,\frac{1}{3}] \cup \{\frac{1}{2}\}\), then \(\frac{1}{\lambda} - 1 \geq 1\) so that \(h^{\frac{1}{\lambda} - 1} \in C_W^{0,\lambda}([0,\delta]) \). Finally, since \(\xi \) and \(\sin \frac{v}{2}\) are continuously differentiable, it follows from above that \(u_W(\cdot, \tau_0) \in C_W^{0,\lambda}([0,\delta]) \) and \(x_Z(\cdot, \tau_0) \in C_Z^{0, \frac{\lambda}{1+\lambda}}([-\delta,\delta])\).
\qed

Now we can Taylor expand \(u \) and \(x\) in terms of \(W\) and \(Z\). When \(\lambda \in (0,\frac{1}{3}] \cup \{\frac{1}{2}\}\), we have  \(u_\tau \in L^\infty\). Hence, 
\[ 
u(W,\tau) = u(W, \tau_0) + \mathcal{O}(|\tau - \tau_0|).
\] 
Then, we can expand \(u(W,\tau_0) \) at \(W = 0\) to get 
\[ 
	u(W,\tau_0) = u(0,\tau_0) + u_W(0,\tau_0) W + \mathcal{O}(|W|^{1 + \lambda}).
\] 
Thus,
\begin{equation}
u(W,\tau) = u(0,\tau_0) + (-1)^{\frac{1}{\lambda} - 1} \frac{\lambda \xi_0 v^{\frac{1}{\lambda} -1}_Y(Y_0,\tau_0) }{2^{\frac{1}{\lambda} - 1}}W + \mathcal{O}(|W|^{1 + \lambda} + |\tau - \tau_0|).
\end{equation}
We play the same game with \(x(Z,\tau)\). By Taylor expanding at \(\tau = \tau_0\), 
\[ 
	x(Z,\tau) = x(Z,\tau_0) + x_\tau(Z,\tau_0)(\tau-\tau_0) + \mathcal{O}(|\tau - \tau_0|^2)
\] 
then writing 
\[ 
		x(Z,\tau_0) = x_0 + x_Z(0,\tau_0)Z + \mathcal{O}(|Z|^{1 + \frac{\lambda}{1 + \lambda}}).
\] 
Hence 
\begin{equation}
		x(Z,\tau) = x_0 + (-1)^{\frac{1}{\lambda}} \frac{\lambda \xi_0 v^{\frac{1}{\lambda}}_Y(Y_0,\tau_0)}{2^{\frac{1}{\lambda}}(1+\lambda)}Z + x_\tau(Z,\tau_0)(\tau-\tau_0) + \mathcal{O}(|Z|^{1 + \frac{\lambda}{1+\lambda}} + |\tau-\tau_0|^2).
\end{equation}
We can now go back to \(Y\) to write  
\begin{equation}\label{u_taylor_I}
		u(Y,\tau) = u(Y_0,\tau_0) + (-1)^{\frac{1}{\lambda} -1} \frac{\lambda \xi_0 v_Y^{\frac{1}{\lambda } - 1}(Y_0,\tau_0) }{2^{\frac{1}{\lambda } -1}}(Y-Y_0)^{\frac{1}{\lambda}} + \mathcal{O}(|Y-Y_0|^{\frac{1+\lambda}{\lambda}} + |\tau - \tau_0|)
\end{equation}
and 
\begin{multline}
x(Y,\tau) = x_0 + (-1)^{\frac{1}{\lambda}} \frac{\lambda \xi_0 v_Y^{\frac{1}{\lambda}}(Y_0,\tau_0)}{2^{\frac{1}{\lambda}}(1+\lambda)}(Y-Y_0)^{\frac{1+\lambda}{\lambda}}  + x_\tau(Y,\tau_0)(\tau-\tau_0)\\ + \mathcal{O}(|Y-Y_0|^{\frac{1+2\lambda}{\lambda}} + |\tau-\tau_0|^2).
\end{multline}
Next, we work on the Taylor expansion \(x_\tau(Y,\tau_0)\) in terms of \(W\). We know \(x_\tau(0,\tau_0) = f'(u(0,\tau_0)) \) and can compute
\[ 
		x_{W\tau}(W,\tau_0) = \frac{\lambda}{2} f''(u(W,\tau_0)) \xi(W,\tau_0) \sin \frac{v(W,\tau)}{2} \left(\frac{\cos \frac{v(W,\tau_0)  }{2}}{W^{\lambda}} \right)^{\frac{1}{\lambda} - 1}.
\] 
Hence 
\[ 
		x_{W\tau}(0,\tau_0):= \lim_{W \to 0^+} x_{W\tau}(W,\tau_0) =  \frac{\lambda \xi_0 f''(u(Y_0,\tau_0))(-v_{Y}(Y_0,\tau_0))^{\frac{1}{\lambda}-1} }{2^{\frac{1}{\lambda}}(1+\lambda) }
\] 
is finite. This implies 
\[ 
	x_\tau(W,\tau_0) = f'(u(0,\tau_0)) + \mathcal{O}(|W|), \quad x_\tau(Y,\tau_0) = f'(u(Y_0,\tau_0)) + \mathcal{O}(|Y-Y_0|^{\frac{1}{\lambda}}). 
\] 
Therefore,
\begin{multline}
		x(Y,\tau) = x_0 +   f'(u(Y_0,\tau_0))(\tau-\tau_0)   + (-1)^{\frac{1}{\lambda}} \frac{\lambda \xi_0 v^{\frac{1}{\lambda}}_Y(Y_0,\tau_0)}{2^{\frac{1}{\lambda}}(1+\lambda)}(Y-Y_0)^{\frac{1+\lambda}{\lambda}}  \\ + \mathcal{O}(|Y-Y_0|^{\frac{1+2\lambda}{\lambda}} + |Y-Y_0|^{\frac{1}{\lambda}}|\tau-\tau_0| +    |\tau-\tau_0|^2).
\end{multline}
We can now solve for \(Y-Y_0\): 
\begin{multline}\label{Y-Y_0-I}
		Y-Y_0 = (-1)^{\frac{1}{1+\lambda}} \left( \frac{2^{\frac{1}{\lambda}}(1+\lambda)  }{\lambda\xi_0 v^{\frac{1}{\lambda}}_Y(Y_0,\tau_0)}\right)^{\frac{\lambda}{1+\lambda}}[x-x_0 - f'(u(Y_0,\tau_0))(\tau-\tau_0)]^{\frac{\lambda}{1+\lambda}} \\
		+ \mathcal{O}(|Y-Y_0|^{\frac{1+2\lambda}{1 + \lambda}} + |Y-Y_0|^{\frac{1}{1+\lambda}}|\tau-\tau_0|^{\frac{\lambda}{1+\lambda}} +    |\tau-\tau_0|^{\frac{2\lambda}{1+\lambda}}).
\end{multline}
Finally, substituting \eqref{Y-Y_0-I} into \eqref{u_taylor_I} gives 
\begin{multline}\label{u_I_+}
		u(x,t) = u(x_0,t_0) - \left(\frac{2\lambda \xi_0 (1+\lambda)^{\frac{1}{\lambda}} }{v_{Y}(Y_0,\tau_0)} \right)^{\frac{\lambda}{1+\lambda}}[x - x_0 - f'(u(x_0,t_0))(t-t_0)]^{\frac{1}{1+ \lambda}} \\
		+ \mathcal{O}(|x-x_0| + |t-t_0|     ).
\end{multline}
This gives our desired expansion in Theorem \ref{main_thm} for even \(n\). Note that for \(Y-Y_0 < 0\), we instead have 
\[
Y_W = -\frac{\lambda}{W^{1-\lambda}}. 
\] 
But by following the same computations as above, we see that
\begin{equation}\label{u_I_neg} 
	u(Y,\tau) = u(Y_0,\tau_0) -  \frac{\lambda \xi_0(v_Y(Y_0,\tau_0))^{\frac{1}{\lambda} -1}}{2^{\frac{1}{\lambda } -1}}(Y-Y_0)^{\frac{1}{\lambda}} + \mathcal{O}(|Y-Y_0|^{\frac{1+\lambda}{\lambda}} + |\tau - \tau_0|). 
\end{equation}
In this case, subbing \eqref{Y-Y_0-I} into \eqref{u_I_neg} still yields \eqref{u_I_+}. 

Now suppose that  \(n \) is odd, so that \(\frac{1}{\lambda} \) is odd and \(\frac{1+\lambda}{\lambda}\) is even. Direct computation still shows that 
\[ 
		u(Y,\tau) = u(Y_0,\tau_0) + (-1)^{\frac{1}{\lambda} -1} \frac{\lambda \xi_0 v_Y^{\frac{1}{\lambda } - 1}(Y_0,\tau_0) }{2^{\frac{1}{\lambda } -1}}(Y-Y_0)^{\frac{1}{\lambda}} + \mathcal{O}(|Y-Y_0|^{\frac{1+\lambda}{\lambda}} + |\tau - \tau_0|)
\] 
and 
\begin{multline*} 
		x(Y,\tau) = x_0 +   f'(u(Y_0,\tau_0))(\tau-\tau_0)   + (-1)^{\frac{1}{\lambda}} \frac{\lambda \xi_0 v^{\frac{1}{\lambda}}_Y(Y_0,\tau_0)}{2^{\frac{1}{\lambda}}(1+\lambda)}(Y-Y_0)^{\frac{1+\lambda}{\lambda}}  \\ + \mathcal{O}(|Y-Y_0|^{\frac{1+2\lambda}{\lambda}} + |Y-Y_0|^{\frac{1}{\lambda}}|\tau-\tau_0| +    |\tau-\tau_0|^2).
\end{multline*}
However, since \(\frac{1+\lambda}{\lambda}\) is even, we only have
\begin{multline} \label{abs_Y}
		|Y-Y_0| =   \left( \frac{2^{\frac{1}{\lambda}}(1+\lambda)  }{\lambda\xi_0 |v_Y(Y_0,\tau_0)|^{\frac{1}{\lambda}}}\right)^{\frac{\lambda}{1+\lambda}}|x-x_0 - f'(u(Y_0,\tau_0))(\tau-\tau_0)|^{\frac{\lambda}{1+\lambda}} \\
		+ \mathcal{O}(|Y-Y_0|^{\frac{1+2\lambda}{1 + \lambda}} + |Y-Y_0|^{\frac{1}{1+\lambda}}|\tau-\tau_0|^{\frac{\lambda}{1+\lambda}} +    |\tau-\tau_0|^{\frac{2\lambda}{1+\lambda}}).
\end{multline} 
Hence, rearranging \eqref{u_taylor_I} and subbing \eqref{abs_Y} gives 
\begin{multline}
		|u(x,t) - u(x_0,t_0)| =  \left(\frac{2\lambda \xi_0 (1+\lambda)^{\frac{1}{\lambda}} }{|v_{Y}(Y_0,\tau_0)|} \right)^{\frac{\lambda}{1+\lambda}}|x - x_0 - f'(u(x_0,t_0))(t-t_0)|^{\frac{1}{1+ \lambda}} + \mathcal{O}(|x-x_0| + |t-t_0|     ).
\end{multline}
This completes the proof of Theorem \eqref{main_thm} for the Type I case. 

\subsection{Type II singularity } 
We can run the same approach for Type II singularities. This time, we introduce the variables 
\begin{equation} 
		W = (Y-Y_0)^{\frac{2-\lambda}{\lambda}}, \quad Z = (Y-Y_0)^{\frac{2+\lambda}{\lambda}}.
\end{equation} 
Note that when \(\lambda = \frac{1}{n}\), we have \(\frac{2- \lambda}{\lambda} = 2n-1 \) and \(\frac{2+\lambda}{\lambda} = 2n + 1\), both of which are odd. Hence, 
\begin{equation} 
		Y_{W} = \frac{\lambda}{(2-\lambda)W^{1 - \frac{\lambda}{2-\lambda}}   }, \quad Y_{Z} = \frac{\lambda}{(2+\lambda)Z^{1 - \frac{\lambda}{2+\lambda}}}.
\end{equation} 
away from \(W = 0\) and \(Z = 0\). We can use the  chain rule to compute 
\begin{align*}
		u_W(W,\tau) & =  \frac{\lambda}{2-\lambda}\xi(W,\tau)  \sin \frac{v(W ,\tau)}{2} \frac{(\cos \frac{v(W, \tau)}{2})^{\frac{1}{\lambda} - 1}}{W^{\frac{2(1-\lambda)}{2-\lambda}}}\\
				 & =  \frac{\lambda}{2-\lambda} \xi(W, \tau) \sin \frac{v(W, \tau)}{2} \left(\frac{\cos \frac{v(W ,\tau)}{2})}{W^{\frac{2 \lambda }{2-\lambda}}}\right)^{\frac{1}{\lambda} - 1}
\end{align*}
and 
\begin{align*}
x_Z(Z,\tau) & = \frac{\lambda}{2+\lambda}\xi(Z,\tau)  \frac{(\cos \frac{v(Z ,\tau)  }{2})^{\frac{1}{\lambda}}    }{Z^{\frac{2}{2+\lambda}}}\\
		 & =\frac{\lambda}{2+\lambda} \xi(Z,\tau) \left(\frac{\cos \frac{v(Z ,\tau)  }{2}}{ Z^{\frac{2\lambda}{2+\lambda}}} \right)^{\frac{1}{\lambda}}.
\end{align*}	
At \(\tau=\tau_0\), we can use L'Hopital's rule to compute 
\begin{align*}
		u_W(0,\tau_0) & := \lim_{W \to 0 } u_W(W,\tau_0) \\
				   & = \frac{\lambda \xi_0}{2-\lambda}\left(\lim_{W \to 0}  \frac{\cos \frac{v(W,\tau_0)}{2})}{W^{\frac{2 \lambda }{2-\lambda}}} \right)^{\frac{1}{\lambda} - 1}\\
				   & = \frac{\lambda \xi_0}{2-\lambda} \left(\lim_{W \to 0} \frac{-\frac{\lambda}{2(2-\lambda)}(\sin \frac{v(W,\tau_0)}{2}) v_Y(W ,\tau_0)  W^{\frac{\lambda}{2-\lambda} - 1}  }{\frac{2\lambda}{2-\lambda} W^{\frac{2\lambda}{2-\lambda}-1}  } \right)^{\frac{1}{\lambda} - 1}\\
				   & = \frac{\lambda \xi_0}{2-\lambda} \left(-\frac{1}{4} \lim_{W \to 0} \frac{ v_Y(W ,\tau_0)  }{ W^{\frac{\lambda}{2-\lambda}}  } \right)^{\frac{1}{\lambda} - 1}\\
				   & = \frac{\lambda \xi_0}{2-\lambda} \left(-\frac{1}{4} \lim_{W \to 0} \frac{ \frac{\lambda}{2-\lambda}v_{YY}(W ,\tau_0)W^{\frac{\lambda}{2-\lambda}-1}   }{ \frac{\lambda}{2-\lambda} W^{\frac{\lambda}{2-\lambda} - 1}  } \right)^{\frac{1}{\lambda} - 1}\\
				   & = (-1)^{\frac{1}{\lambda} -1} \frac{\lambda \xi_0}{4^{\frac{1}{\lambda} - 1}(2-\lambda) }v^{\frac{1}{\lambda} - 1}_{YY}(Y_0,\tau_0) \neq 0
\end{align*}
and 
\begin{align*}
		x_Z(0,\tau_0) & := \lim_{Z \to 0 } x_Z(Z,\tau_0) \\
				   & = \frac{\lambda \xi_0}{2+\lambda} \left(\lim_{Z \to 0} \frac{\cos \frac{v(Z ,\tau_0)  }{2}}{Z^{\frac{2\lambda}{2+\lambda}} } \right)^{\frac{1}{\lambda}}\\
				   & = (-1)^{\frac{1}{\lambda}} \frac{\lambda \xi_0}{4^{\frac{1}{\lambda}}(2+\lambda)}v_{YY}^{\frac{1}{\lambda}} (Y_0,\tau_0)  \neq 0.
\end{align*}
We have a similar proposition to the Type I case. 
\begin{proposition}\label{holder_II}
		Suppose  \(\lambda = \frac{1}{n}\) where \(n\geq 2\) is an integer.  Then there is some \(\delta > 0\) such that \(u_W(\cdot,\tau_0) \in C_W^{0, \frac{\lambda}{2-\lambda} }([-\delta,\delta]) \) and \(x_Z(\cdot,\tau_0) \in C_Z^{0,\frac{\lambda}{2+\lambda}}([-\delta, \delta])\). 
\end{proposition}	
\pf
	The proof follows similarly as Proposition \ref{holder_I}. Here we define
	\[ 
			h(Y) := \frac{\cos \frac{v(Y,\tau_0)}{2}}{(Y-Y_0)^2}
	\] 
	since 
	\[ 
			u_W(W,\tau_0) = \frac{\lambda}{2-\lambda} \xi(W,\tau_0) \sin \frac{v(W,\tau_0)}{2} (h(W^{\frac{2\lambda}{2-\lambda}} + Y_0))^{\frac{1}{\lambda} -1} 
	\] 
	and 
	\[ 
			x_Z(Z,\tau_0) = \frac{\lambda}{2+\lambda} \xi(Z,\tau_0) \left(h(Z^{\frac{2\lambda}{2+\lambda}} + Y_0) \right)^{\frac{1}{\lambda}}.
	\] 
	Direct computation shows 
	\[ 
			\lim_{Y \to Y_0} h'(Y)  = -\frac{1}{12}v_{YYY}(Y_0,\tau_0).
	\] 
	Hence, \(h\) is Lipschitz near \(Y = Y_0\) and the result follows. 
\qed

As a result, we can get the expansions
\begin{equation} 
		u(W,\tau) = u(0,\tau_0) + (-1)^{\frac{1}{\lambda} - 1} \frac{\lambda \xi_0 v_{YY}^{\frac{1}{\lambda} - 1}(Y_0,\tau_0) }{4^{\frac{1}{\lambda} - 1}(2-\lambda) }W  + \mathcal{O}(|W|^{1 + \frac{\lambda}{2-\lambda}} + |\tau-\tau_0| )
\end{equation} 
and 
\begin{align}
		x(Z,\tau) = x_0 + (-1)^{\frac{1}{\lambda}} \frac{\lambda \xi_0 v^{\frac{1}{\lambda}}_{YY}(Y_0,\tau_0) }{4^{\frac{1}{\lambda}}(2+\lambda)}Z + x_\tau(Z,\tau_0)(\tau - \tau_0)  + \mathcal{O}(|Z|^{1 + \frac{\lambda}{2+\lambda}} + |\tau-\tau_0|^2).
\end{align}
Thus
\begin{align}\label{u+} 
		u(Y,\tau) = u(Y_0,\tau_0)  + (-1)^{\frac{1}{\lambda} - 1} \frac{\lambda \xi_0 v^{\frac{1}{\lambda} - 1}_{YY}(Y_0,\tau_0)  }{4^{\frac{1}{\lambda} - 1}(2-\lambda) }(Y-Y_0)^{\frac{2-\lambda}{\lambda}}  + \mathcal{O}(|Y-Y_0|^{\frac{2}{\lambda}}  + |\tau-\tau_0| )
\end{align} 
and 
\begin{multline}
		x(Y,\tau) = x_0 + (-1)^{\frac{1}{\lambda}} \frac{\lambda \xi_0 v^{\frac{1}{\lambda}}_{YY}(Y_0,\tau_0) }{4^{\frac{1}{\lambda}}(2+\lambda)}(Y-Y_0)^{\frac{2+\lambda}{\lambda}} + x_\tau(Y,\tau_0)(\tau-\tau_0) \\+ \mathcal{O}(|Y-Y_0|^{\frac{2+2\lambda}{\lambda}} + |\tau-\tau_0|^2). 
\end{multline}
To investigate the term \(x_\tau(Y,\tau_0) \) in  \(W\), we calculate
\begin{align*} 
		x_{W\tau}(W,\tau_0)  & = \frac{\lambda}{2(2-\lambda)} f''(u(W,\tau_0)) \xi(W,\tau_0) \sin \frac{v}{2} \left(\frac{\cos \frac{v(W ,\tau_0)}{2})}{W^{\frac{2 \lambda }{2-\lambda}}}\right)^{\frac{1}{\lambda} - 1}
\end{align*} 
with 
\[ 
		x_{W\tau}(0,\tau_0)  := \lim_{W \to 0}x_{W\tau}(W,\tau_0) = \frac{\lambda \xi_0 f''(u(Y_0,\tau_0))(-v_{YY}(Y_0,\tau_0))^{\frac{1}{\lambda} -1}}{4^{\frac{1}{\lambda}-1}2(2-\lambda)} < \infty. 
\] 
Hence 
\[ 
		x_\tau(W,\tau_0) = f'(u(0,\tau_0)) + \mathcal{O}(|W|), \quad x_\tau(Y,\tau_0) = f'(u(Y_0,\tau_0)) + \mathcal{O}(|Y-Y_0|^{\frac{2-\lambda}{\lambda}})
\] 
and we can write 
\begin{multline} 
		x(Y,\tau) = x_0 + f'(u(Y_0,\tau_0))(\tau-\tau_0) + (-1)^{\frac{1}{\lambda}}\frac{\lambda \xi_0 v_{YY}^{\frac{1}{\lambda}}(Y_0,\tau_0) }{4^{\frac{1}{\lambda}}(2+\lambda)}(Y-Y_0)^{\frac{2+\lambda}{\lambda}} \\ 
	+ \mathcal{O}(|Y-Y_0|^{\frac{2+2\lambda}{\lambda}} + |\tau-\tau_0|^2 + |Y-Y_0|^{\frac{2-\lambda}{\lambda}}|\tau-\tau_0| ).
\end{multline} 
Now we solve for \(Y-Y_0\): 
\begin{multline}\label{Y+}
		Y-Y_0  = (-1)^{\frac{1}{2+\lambda}} \left( \frac{4^{\frac{1}{\lambda}}(2+\lambda)}{\lambda \xi_0 v^{\frac{1}{\lambda}}_{YY}(Y_0,\tau_0)} \right)^{\frac{\lambda}{2+\lambda}} [x - x_0 -f'(u(Y_0,\tau_0))(\tau-\tau_0)]^{\frac{\lambda}{2+\lambda}}\\
		+ \mathcal{O}(|Y-Y_0|^{\frac{2+2\lambda}{\lambda}} + |\tau-\tau_0|^{\frac{2\lambda}{2+\lambda}} + |Y-Y_0|^{\frac{2-\lambda}{2+\lambda}}|\tau-\tau_0|^{\frac{\lambda}{2+\lambda}}  ).
\end{multline}
Finally, substituting \eqref{Y+} into \eqref{u+} yields
\begin{multline}\label{u1}
		u(x,t) = u(x_0,t_0) -  \left(\frac{4\lambda^2 \xi_0^2 (2+\lambda)^{\frac{2-\lambda}{\lambda}}  }{v_{YY}(Y_0,\tau_0)(2-\lambda)^{\frac{2+\lambda}{\lambda}}    } \right)^{\frac{\lambda}{2+\lambda}}[x-x_0 - f'(u(x_0,t_0))(t-t_0)]^{\frac{2-\lambda}{2+\lambda}} \\
		+ \mathcal{O}((|x-x_0| + |t-t_0|)^{\frac{2}{2+\lambda} } ).
\end{multline} 
This completes the proof of Theorem \ref{main_thm} when \(g(u,u_x) = 0\). 

\subsection{Hölder Continuity}

Theorem \ref{main_thm} shows that generic regular solutions are Hölder continuous of exponent \(\frac{2-\lambda}{2+\lambda}\) near the singular curves when \(\lambda = \frac{1}{n}\) for \(n \geq 2\) an integer. This matches expectations for the regularity previously studied for Camassa-Holm and other nonlinear wave equations in \cite{BressanHuangYu,LZ,KKY}. It is reasonable to suspect that this regularity holds for all \(\lambda \in (0,1)\). In fact, the method above also works for other \(\lambda\), such as when \(\lambda = \frac{m}{n}\) where \(m\) is odd and \(\lambda\) is small enough so that \(v_{YYY}\) is bounded. For other \(\lambda\), such as when \(\lambda\) is irrational, more care is needed to handle the change of variables when \(Y-Y_0 < 0\).

\section{General Camassa-Holm type equation with $g(u,u_x)\neq 0$} 

When $\lambda=\frac{1}{2a}$, $f'(u)=u^a$ in \eqref{spe} with any positive integer $a$ and $g(u, u_x)=u^{a+1}-P-Q_x$:
\beq \label{E-2-2}
\begin{array}{lll}u_{t x}+u^a u_{xx}+\frac{1}{2} u^{a-1} u_x^2- u^{a+1}+P+Q_x=0.
 \end{array}
 \eeq
where
 $P:=p*[\frac{2a-1}{2}u^{a-1} u_x^2+ u^{a+1}] \ \ {\rm and}\ \ Q:=\frac{a -1}{2}\,p*[ u^{a-2}u_x^3].$
 This equation gives Camassa-Holm equation when $a=1$ and Novikov equation when $a=2$. We will consider the Cauchy problem for initial data \(u_0(x) \in H^1(\mathbb{R}) \cap W^{1,\frac{1}{\lambda}}(\mathbb{R}) \). 

 \subsection{Semi-linear System} 
 The equation of the characteristic of \eqref{E-2-2} is
 \[ 
		\frac{dx^c(s)}{ds} = u^a(x^c(s),s)). 
 \] 
 We define the change of coordinates \((x,t) \to (Y,\tau)\) as 
 \[ 
		Y = \int_0^{x^c(0,x,t)}(1+u_x^2(y,0))^{\frac{1}{2\lambda}} dy, \quad \tau = t. 
 \] 
 Similar as above, we compute 
\begin{equation} \label{u_Y_nz} 
		u_Y  = \frac{u_x}{Y_x} = \xi \frac{u_x}{(1+u_x^2)^a } = \frac{1}{2}\xi \sin v \cos^{2 a - 2} \frac{v}{2}  
\end{equation}
and 
\[ 
		u_\tau = u_t + u^a u_x = -P_x - Q. 
\] 
Moreover, 
\begin{align*}
	v_\tau & = \frac{2}{1+u_x^2}(u_{xt} + u^a u_{xx}) \\
		   & = -u^{a-1} \frac{u_x^2}{1+u_x^2} + 2u^{a+1} \frac{1}{1+u_x^2} - \frac{2}{1+u_x^2} (P+Q_x) \\
		   & = - u^{a-1} \sin^2\frac{v}{2} + 2u^{a+1} \cos^2 \frac{v}{2} - 2\cos^2 \frac{v}{2}(P+Q_x) 
\end{align*}
and 
\begin{align*}
	\xi_\tau & = \frac{1}{Y_x}a(1+u_x^2)^{a-1} (u_x^2)_\tau - \frac{(1+u_x^2)^a}{Y_x^2}(Y_x)_\tau \\
			 & = \frac{a }{Y_x}(1+u_x^2)^{a- 1}( (u_x^2)_t + u^a (u_x^2)_x + \frac{1}{a} (1+u_x^2) au^{a-1}u_x     )  \\
			 & = \frac{(1+u_x^2)^a }{Y_x} \frac{a u^{a-1}u_x}{1+u_x^2} + \frac{(1+u_x^2)^a}{Y_x}\frac{u_x}{1+u_x^2} \frac{1}{u_x^{2a-1}}( (u_x^{2a})_t + (u^a u_x^{2a})_x ) \\
			 & = \frac{a}{2} u^{a-1} \xi \sin v + a u^{a+1} \xi \sin v - a \xi \sin v(P+Q_x).
\end{align*}
Our semilinear system becomes
\begin{equation} \label{nlsemeq2}
	\begin{cases}
 u_\tau & = -P_x - Q, \\
 v_\tau & = - u^{a-1} \sin^2\frac{v}{2} + 2u^{a+1} \cos^2 \frac{v}{2} - 2\cos^2 \frac{v}{2}(P+Q_x) ,\\
 \xi_\tau & = \frac{a}{2} u^{a-1} \xi \sin v + a u^{a+1} \xi \sin v - a \xi \sin v(P+Q_x).
	\end{cases}
\end{equation} 
with the initial data 
\begin{equation}\label{id_nonzero}
	\begin{cases}
		u(Y,0) := u_0(x(Y,0)), \\
		v(Y,0) := 2\arctan (u'_0(x(Y,0))),\\
		\xi(Y,0) := 1. 
	\end{cases}
\end{equation}
Similarly,
\begin{equation}\label{xtyn}
		\begin{cases}
			x_Y & =  \xi \cos^{2a}\frac{v}{2}, \\
			x_\tau & = u^a,
			\end{cases} \quad \begin{cases}
			t_Y & = 0, \\
			t_\tau & = 1. 
	\end{cases}
\end{equation} 
Using this seminar system, we can prove Theorem \ref{gr_nonzero}. We omit proof as it follows in the same way as the proof of Theorem \ref{gr_zero}.  For the generic regulairty proof, it is critical that \(P,P_x,Q,Q_x\) are Lipschitz. This has been demonstrated in \cite{CSZ} when \(a = 1\) or \(a \geq 2\). Moreover, we can see from \eqref{u_Y_nz} and \eqref{xtyn} that the expansion near the singular curve follows in exactly same as the case where \(g(u,u_x) = 0\). Thus, the second part of Theorem \eqref{main_thm} follows immediately, and the constants \(C_1(\lambda)\) and \(C_2(\lambda)\) remain the same. For demonstration, we will investigate the singular expansions of the Camassa-Holm equation (\(a=1\)) and Novikov equation (\(a=2\)). Note that \(\tan\frac{v}{2} = u_x\), so we have a singular at \(2\pi\) multiples of \(v = \pi\). For Type II singularities, we assume that the point \( (Y_0,\tau_0)\)
\beq
v(Y_0,\tau_0) = \pi,\quad v_Y(Y_0,\tau_0) = 0,\quad v_{YY}(Y_0,\tau_0) \neq 0.
\eeq
We will only look at the expansion near the Type II singularities to derive the order of the cusp singularities. 

\subsection{Camassa-Holm: \( a = 1\)}
When \(a = 1\),  we have 
\[ 
	\begin{cases}
		u_Y = \frac{1}{2}\xi \sin v, \\
		u_\tau = -P_x,  \\
		x_Y = \xi \cos^2 \frac{v}{2}, \\
		x_\tau = u. 
	\end{cases}
\] 
The computations for the singular behavior about the Type I and Type II singularities will match the  general Hunter-Saxton case above where \(\lambda = \frac{1}{2}\). In fact, we can also compute the mixed term \(u_{Y\tau}\) for more detail. Using \(P-P_{xx} = \frac{1}{2}u_x^2 - u^2\), 
\begin{align*} 
		u_{Y\tau} & = (-P_x)_Y = \frac{-P_{xx}}{Y_x} = -\xi (\cos^{2}\frac{v}{2})(P - \frac{1}{2}u_x^2 - u^2)= \frac{\xi}{2} \sin^2 \frac{v}{2}  + \xi \cos^2\frac{v}{2}(u^2 - P)
\end{align*} 
So \(u_{Y\tau}(Y_0,\tau_0) = \frac{1}{2}\xi_0 \). For Type II singularities, it follows from \eqref{u+} that  
\begin{equation}\notag
\begin{split}
		u(Y,\tau) = u(Y_0,\tau_0)  -P_x(Y_0,\tau_0)(\tau-\tau_0) - \frac{v_{YY}(Y_0,\tau_0)\xi_0}{12} (Y-Y_0)^3 
+ \frac{1}{2}\xi_0(Y-Y_0)(\tau-\tau_0)  \\ 
		+ \mathcal{O}(|Y-Y_0|^4 + |\tau-\tau_0|^2).
\end{split}
\end{equation}
Thus
\begin{equation}\notag
\begin{split}
				u(x,t) = u(x_0,t_0) - \frac{1}{12} \left(\frac{80^3\xi_0^2 }{v_{YY}(Y_0,\tau_0)}    \right)^{1/5}[(x-x_0) - u(x_0,t_0)(t-t_0)]^{3/5} \\ 
				+ \mathcal{O} \left(|t-t_0| + |x-x_0|\right)^{4/5}.
\end{split}
\end{equation}
Here we see the solution is Hölder continuous of exponent 3/5 near the singular curve. Our expansion and order of singularity agree with the previous work of  \cite{LZ, KKY, BressanHuangYu}.

\subsection{Novikov: \(a = 2\) }
When \(a = 2\), we have 
\[ 
	\begin{cases}
			u_Y = \frac{1}{2}\xi \sin v \cos^{2}\frac{v}{2},   \\
		u_\tau = -P_x - Q, \\ 
		x_Y = \xi \cos^{4}\frac{v}{2}, \\
		x_\tau = u^2.
	\end{cases}	
\] 
We can also compute the Taylor expansion in the same was as we did before with  \(\lambda = \frac{1}{4}\). Again, we will consider the Type II case and can compute the mixed partials for more detail. Here 
\begin{align*}
		u_{Y\tau} & = -\xi \cos^4\frac{v}{2}(P_{xx} +  Q_x ) \\ 
			   & = -\xi \cos^4\frac{v}{2}(P - \frac{3}{2}u u_x^2 - u^3 + Q_x) \\
			   & = -\xi \cos^4 \frac{v}{2}(P + Q_x) + \frac{3}{2}\xi u \sin^2 \frac{v}{2} \cos^2 \frac{v}{2}  + \xi u^3 \cos^4 \frac{v}{2}  \label{uyt} 
\end{align*}
At \( (Y_0,\tau_0)\) we have \(u_{Y\tau}(Y_0,\tau_0) = 0\). Now 
\begin{equation}\notag
\begin{split}
		 u_{YY\tau}  = &\frac{3}{2}\left( \xi_Y u \sin^2\frac{v}{2} \cos^2\frac{v}{2} + \xi u_Y \sin^2\frac{v}{2} \cos^2\frac{v}{2}  + \xi u v_Y \sin\frac{v}{2} \cos^3 \frac{v}{2} 
				 - \xi u v_Y \sin^3\frac{v}{2}\cos\frac{v}{2}\right) \\
				& + \xi_Y u^3 \cos^4 \frac{v}{2} + 3\xi u^2 u_Y \cos^4 \frac{v}{2} - 2\xi u^3v_Y \cos^3 \frac{v}{2}\sin\frac{v}{2}\\
				& - \xi_Y \cos^4\frac{v}{2}(P+Q_x) + 2\xi v_Y\cos^3\frac{v}{2}\sin\frac{v}{2}(P+Q_x) - \xi^2 \cos^{8}\frac{v}{2} (P_x + Q_{xx}).
\end{split}
\end{equation}
But \(Q_{xx} = Q - \frac{1}{2}u_x^3 = Q - \frac{1}{2} \tan^3\frac{v}{2}\). Hence 
\begin{equation}
\begin{split}
		 u_{YY\tau}  = &\frac{3}{2}\left( \xi_Y u \sin^2\frac{v}{2} \cos^2\frac{v}{2} + \xi u_Y \sin^2\frac{v}{2} \cos^2\frac{v}{2}  + \xi u v_Y \sin\frac{v}{2} \cos^3 \frac{v}{2}  - \xi u v_Y \sin^3\frac{v}{2}\cos\frac{v}{2}\right) \\
				& + \xi_Y u^3 \cos^4 \frac{v}{2} + 3\xi u^2 u_Y \cos^4 \frac{v}{2} - 2\xi u^3v_Y \cos^3 \frac{v}{2}\sin\frac{v}{2}\\
				& - \xi_Y \cos^4\frac{v}{2}(P+Q_x) + 2\xi v_Y\cos^3\frac{v}{2}\sin\frac{v}{2}(P+Q_x) - \xi^2 \cos^{8}\frac{v}{2} (P_x + Q) \\
				&  + \frac{1}{2}\xi^2 \cos^5 \frac{v}{2} \sin^3\frac{v}{2}.   \label{uyyt} 
\end{split}
\end{equation}
Again we see that \(u_{YY\tau}(Y_0,\tau_0) = 0\). Because of the \(- \frac{3}{2} \xi u v_Y \sin^3 \frac{v}{2} \cos \frac{v}{2}\) term, repeated differentiation shows that 
\[ 
	 u_{YYY\tau}(Y_0,\tau_0) = 0, \quad  u_{YYYY\tau}(Y_0,\tau_0) = 0, \quad  u_{YYYYY\tau}(Y_0,\tau_0) = \frac{3}{2}\xi_0 u(x_0,t_0) v_{YY}^2(Y_0,\tau_0) \neq 0.  
\]

By \eqref{u+}, 
\begin{align*}
		 u(Y,\tau)  =& u(Y_0,\tau_0) - ( Q(Y_0,\tau_0) + P_x(Y_0,\tau_0))(\tau-\tau_0)  -\frac{\xi_0v_{YY}^{3}(Y_0,\tau_0)}{448}(Y-Y_0)^{7}\\
		& + \frac{45}{1440}\xi_0 u(x_0,t_0) v_{YY}^2(Y_0,\tau_0)(Y-Y_0)^5(\tau-\tau_0) + \mathcal{O}(|Y-Y_0|^{8} + |\tau-\tau_0|^2).
\end{align*}
Therefore 
\begin{align*}
u(x,t) =& u(x_0,t_0)  - \left( \frac{ 4 ( 9^7)  \xi_0^2 }{ 7^9 v_{YY}(Y_0,\tau_0)}\right)^{\frac{1}{9}}[(x-x_0) - u(x_0,t_0)^2 (t-t_0)]^{\frac{7}{9}}  + \mathcal{O}( (|x-x_0| + |t-t_0|)^{\frac{8}{9}}). 
\end{align*}
In particular, the solution is Hölder continuous of exponent \(7/9\) near the singular curve.

\section*{Conflict of interest statement}
There is no conflict of interest.
\section*{Data availability statement}
No data are used in this paper.

\section*{Acknowledgments}
The second author is partially supported by National Science Foundation  with grant  DMS-2306258 and DMS-2605028.
The fourth author is partially supported by National Science Foundation with grant DMS-2206218.
%

 
\end{document}